\def\DateTime{11/May/1999} % 19:50
\def\Version{Version $1.1$}
\def\yes{\if00}

\def\iftwelvept{\yes}

\iftwelvept
\documentclass[leqno,12pt]{amsart}
\else
\documentclass[leqno]{amsart}
\fi
\usepackage{amssymb}
\usepackage{amscd}
\usepackage{latexsym}

\iftwelvept
\else\fi

\theoremstyle{plain}
\newtheorem{Theorem}{Theorem}[section]
\newtheorem{Proposition}[Theorem]{Proposition}
\newtheorem{Lemma}[Theorem]{Lemma}
\newtheorem{Corollary}[Theorem]{Corollary}
\newtheorem{Claim}{Claim}[Theorem]

\theoremstyle{definition}

\newtheorem{Remark}[Theorem]{Remark}
\newtheorem{Example}[Theorem]{Example}

\renewcommand{\theTheorem}{\arabic{section}.\arabic{Theorem}}

\newcommand{\ZZ}{{\mathbb{Z}}}
\newcommand{\QQ}{{\mathbb{Q}}}
\newcommand{\RR}{{\mathbb{R}}}
\newcommand{\CC}{{\mathbb{C}}}
\newcommand{\PP}{{\mathbb{P}}}
\newcommand{\OO}{{\mathcal{O}}}

\newcommand{\KQ}{\operatorname{[K:\QQ]}}
\newcommand{\achern}{\operatorname{\widehat{c}}}
\newcommand{\chern}{\operatorname{c}}
\newcommand{\adeg}{\operatorname{\widehat{deg}}}
\newcommand{\td}{\operatorname{tr.deg}}
\newcommand{\ord}{\operatorname{ord}}

\newcommand{\Supp}{\operatorname{Supp}}
\newcommand{\Ker}{\operatorname{Ker}}
\newcommand{\Coker}{\operatorname{Coker}}
\newcommand{\Bs}{\operatorname{Bs}}
\newcommand{\Proof}{{\sl Proof.}\quad}
\newcommand{\QED}{{\unskip\nobreak\hfil\penalty50\quad\null\nobreak\hfil
{$\Box$}\parfillskip0pt\finalhyphendemerits0\par\medskip}}

\begin{document}

%%%%%%%%%
% Title %
%%%%%%%%%

\title[periodic points]{A remark on periodic points on varieties \\ 
over a field of finite type over $\QQ$}

\author{Shu Kawaguchi}
\address{Department of Mathematics, Faculty of Science,
Kyoto University, Kyoto, 606-01, Japan}
\email[Shu Kawaguchi]{kawaguch@kusm.kyoto-u.ac.jp}
\keywords{periodic point, rational point, iteration}
\subjclass{14G05, 11G99}
\date{\DateTime, (\Version)}

\begin{abstract}
Let $M$ be a field of finite type over $\QQ$ and  
$X$ a variety defined over $M$. 
We study when 
the set $\{ P \in X(K) \mid f^{\circ n} (P) = P  
\ \text{for some} \ n \geq 1 \}$ 
is finite for any finite extension fields $K$ of $M$ and 
for any dominant $K$-morphisms $f : X \to X$ with 
$\deg f \geq 2$. 
\end{abstract}

\maketitle

\section*{Introduction}
\renewcommand{\theTheorem}{\Alph{Theorem}}
By a variety, we mean 
an integral separated scheme of finite type 
over a ground field. 
Let $M$ be a field of finite type over $\QQ$ 
and $X$ a variety defined over $M$. 

Let $K$ be a finite extension field of $M$ and 
$f : X \to X$ a dominant morphism defined over $K$. 
We say that a point $P \in X(K)$ is periodic with respect to $f$ if 
there is a positive integer $n$ with $f^{\circ n} (P) = P$. 
Let $X(K)_{per,f}$ be the set of periodic $K$-points 
with respect to $f$.  
We say that $X$ is periodically finite if 
$X(K)_{per,f}$ is a finite set 
for any finite extension fields $K$ of $M$ and 
any dominant $K$-morphisms $f : X \to X$ with $\deg f \geq 2$. 

In this paper, we study when $X$ is periodically finite. 

In order to show the finiteness of $X(K)_{per,f}$, 
we introduce the set of backward $K$-orbits of $f$, 
denoted by $\varprojlim_{f} X(K)$, which is defined by 
\[
\varprojlim_{f} X(K) 
= \{
  (x_n)_{n=0}^{\infty} \in \prod_{n=0}^{\infty} X(K) \mid 
  f(x_{n+1}) = x_n \quad (n \geq 0)
  \}.
\] 
It is easy to see that if $\varprojlim_{f} X(K)$ is a finite set, 
then so is $X(K)_{per,f}$ and $ \# \varprojlim_{f} X(K) = \# X(K)_{per,f}$ 
(cf. Lemma~\ref{lemma:elemetary:property:of:periodic:pts}). 

We obtain the following results.

\begin{Theorem}[cf. Corollary~\ref{cor:finite} %
and \S \textup{\ref{sec:final}}]
Let $X$ be a geometrically irreducible normal projective variety 
defined over a field of finite type over $\QQ$. 
Assume that the Picard number of $X$ is $1$ 
(for example, $X$ is $\PP^n$ 
or a geometrically irreducible normal projective curve). 
Then $X$ is periodically finite. 
\end{Theorem}

We prove this by using  
Northcott's finiteness theorem of height functions. 
More precisely, this result is a corollary of the fact that 
if there is an ample line bundle 
$L$ such that $f^{*}(L) \otimes L^{-1}$ is also ample,   
then $\varprojlim_{f} X(K)$ is a finite set 
(Theorem~\ref{theorem:finite}).

We also show:

\begin{Theorem}[cf. Corollary~\ref{cor:curve} %
and \S \textup{\ref{sec:final}}]
Let $C$ be a curve defined over 
a field of finite type over $\QQ$. 
Then $C$ is periodically finite. 
\end{Theorem}

\begin{Theorem}[cf. Theorem~\ref{thm:abelian:variety} % 
and \S \textup{\ref{sec:final}}]
Let $A$ be an abelian variety defined over 
a field of finite type over $\QQ$. 
Then $A$ is periodically finite if and only if 
$A$ is simple.
\end{Theorem}

For a number field case, we also show:

\begin{Theorem}[cf. Theorem~\ref{theorem:surfaces:with:kx:geq:0}]
Let $X$ be a smooth projective surface 
with the non-negative Kodaira dimension 
such that $X$ is defined over 
a number field (for the case of a field of finite type over 
$\QQ$, see \S 6). 
Then $X$ is {\em not} periodically finite if and only if 
$X$ is one of the following types;
\begin{enumerate}
\renewcommand{\labelenumi}{(\roman{enumi})}
\item
$X$ is an abelian surface which is {\em not} simple, or
\item
$X$ is a hyperelliptic surface. 
\end{enumerate}
\end{Theorem}

In order to clarify the argument, 
$M$ is assumed to be a number field before \S \ref{sec:final}, 
where in \S \ref{sec:final}, we deal with 
a field of finite type over $\QQ$ in general. 

The author expresses his gratitude to Prof. Yoshio Fujimoto 
for telling him their results (\cite{FS}, \cite{FS2}).   
While the author was preparing this paper, he was 
partially supported by JSPS Research Fellowships
for Young Scientists.

\medskip
\section{Quick review of height theory}
\renewcommand{\theTheorem}{\arabic{section}.\arabic{Theorem}}
In this section, 
we recall some properties of height functions. 
We refer to \cite{Sil} for details.  

Let $h : \PP^n(\overline{\QQ}) \to \RR$ be 
the logarithmic height function. 
Namely, for a point $x \in \PP^n(\overline{\QQ})$, 
$h(x)$ is defined by 
\[
h(x) = \frac{1}{\KQ} \sum_{v \in M_K} 
         \log \left( \max_{1 \leq i \leq n} \{ \vert x_i \vert_{v} \}
              \right), 
\] 
where $x = (x_0, x_1, \ldots, x_n) \in \PP^n(K)$ is its coordinate 
over a sufficiently large number field $K$, 
and $M_K$ be the set of all places of $K$. 

Now let $X$ be a projective variety defined over $\overline{\QQ}$, 
$\phi : X \to \PP^n$ a morphism over $\overline{\QQ}$. 
For a point $x \in X(\overline{\QQ})$, 
we define the height of $x$ with respect to $\phi$, 
denoted by $h_{\phi}(x)$, to be $h_{\phi}(x) = h(\phi(x))$. 

Then the following theorem holds. 

\begin{Theorem}[Height Machine]
\label{theorem:height:machine}
For every line bundle $L$ on a projective variety $X$ 
defined over $\overline{\QQ}$, 
there exists a unique function 
$h_L : X(\overline{\QQ}) \to  \RR$ modulo bounded functions 
with the following property; 
\begin{enumerate}
\renewcommand{\labelenumi}{(\roman{enumi})}
\item 
For any two line bundles $L_1, L_2$, 
$h_{L_1 \otimes L_2} = h_{L_1} + h_{L_2} + O(1)$.  
\item 
If $f : X \to Y$ be a morphism of projective varieties 
over $\overline{\QQ}$, 
then $h_{f^*(L)} = f^*(h_L) + O(1)$.
\item 
If $\phi : X \to \PP^n$ a morphism over $\overline{\QQ}$, 
then $h_{\phi^*\left(\OO_{\PP^n}(1)\right)} 
= h_{\phi} + O(1)$. 
\end{enumerate}
\end{Theorem}

We also recall some properties of height functions.  

\begin{Theorem}
\label{theorem:positiveness:of:height}
\begin{enumerate}
\renewcommand{\labelenumi}{(\roman{enumi})}
\item
(positiveness)
If we denote 
$\Supp \left( \Coker (H^0(X,L) \otimes \OO_X) \to L \right)$ 
by $\Bs (L)$, then 
$h_L$ is bounded below on $\left(X \setminus \Bs(L) \right)$. 
\item
(Northcott)
Assume $L$ be ample.  
Then for any $d \geq 1$ and $M \geq 0$, 
\[
\{ x \in X(\overline{\QQ}) 
\mid  h_{L}(x) \leq M, \quad [\QQ(x):\QQ] \leq d \}
\] 
is a finite set. 
\end{enumerate}
\end{Theorem} 

For Theorem~\ref{theorem:height:machine}, 
we refer to \cite[Theorem 3.3]{Sil}. 
For Theorem~\ref{theorem:positiveness:of:height}, 
we refer to \cite[Corollary 3.4 and Proposition 3.5]{Sil}. 
Although in \cite{Sil} Theorem~\ref{theorem:height:machine} (ii) is 
written for a morphism of smooth projective varieties, 
it also holds for not necessarily smooth projective varieties. 

\medskip
\section{Finiteness}
Let $X$ be a variety defined over a number field $M$. 
Let $K$ be a finite extension of $M$ and 
$f : X \to X$ a dominant morphism defined over $K$. 

We say that a point $P \in X(K)$ is {\em periodic} with respect to $f$ if 
there is a positive integer $n$ with $f^{\circ n} (P) = P$. 
Let $X(K)_{per,f}$ be the set of periodic $K$-points with respect to $f$.  

We also define the set of {\em backward $K$-orbits} of $f$, 
denoted by $\varprojlim_{f} X(K)$, to be 
\[
\varprojlim_{f} X(K) 
= \{
  (x_n)_{n=0}^{\infty} \in \prod_{n=0}^{\infty} X(K) \mid 
  f(x_{n+1}) = x_n \quad (n \geq 0)
  \}.
\] 

We say that $X$ is {\em periodically finite} if 
for any finite extension $K$ of $M$ and 
for any dominant $K$-morphism $f : X \to X$ 
with $\deg f \geq 2$, $X(K)_{per,f}$ is a finite set. 

In this paper, 
we would like to study what kind of $X$ is periodically finite. 

We first remark elementary properties of 
$X(K)_{per,f}$ and $\varprojlim_{f} X(K)$. 

\begin{Lemma}
\label{lemma:useful}
Let  $S \subset X(K)$ be a finite set 
and $(x_n)_{n=0}^{\infty} \in \varprojlim_{f} X(K)$. 
Assume that there is a subsequence 
$(x_{n_i})_{i=0}^{\infty}$ consisting of elements in $S$. 
Then $(x_n)_{n=0}^{\infty}$ is periodic, i.e., 
there is a positive integer $p$ with 
$x_{n+p} = x_n$ for $n \geq 0$. 
Moreover, $(x_n)_{n=0}^{\infty}$ is 
uniquely determined by $x_0$. 
\end{Lemma}

\Proof
Since $S$ is a finite set, there is an element $s \in X(K)$ 
such that, for infinitely many $n$, $x_n$ equals to $s$. 
Let $(x_{n_j})_{j=0}^{\infty}$ be the subsequence of 
$(x_n)_{n=0}^{\infty}$ with $x_{n_j} = s$ for $j \geq 0$. 
Let us set $p = n_1 -n_0$. We show that $n_2 -n_1 = p$. 
Indeed, 
since $f^{\circ q} (x_{n_2}) = x_{n_1}$, 
if we set $q = n_2 -n_1$, then we have $f^{\circ q}(s) = s$. 
If we assume $q >p$, then $n_2 > n_2 -p > n_1$ and 
$x_{n_2} = x_{n_2 - p} = x_{n_1} = s$. This is a contradiction. 
If we assume $p >q$, then we similarly have a contradiction. 
Thus $n_2 - n_1 = n_1 - n_0 = p$. In the same way, 
$n_{j+1} - n_j = p$ for any $j \geq 0$. 
Now let us take any $n \geq 0$. We fix an $n_j$ with $n_j > n$ and 
set $r = n_j -n$. Then $n_j +p = n_{j+1}$ and 
$n_{j+1} -(n+p) = r$. Therefore, we get 
\[
x_{n+p} = f ^{\circ r}(x_{n_{j+1}}) 
        = f ^{\circ r}(s) = f ^{\circ r}(x_{n_j}) = x_n.   
\]
This shows that $(x_n)_{n=0}^{\infty}$ is periodic. 
Moreover if we divide $n$ by $p$ and write 
$n = q p + l$ with $0 \leq l \leq p-1$, 
then it is easy to see that $x_n = f^{\circ (p-l)}(x_0)$. 
This shows the latter assertion of the lemma.    
\QED

The next lemma gives the relationship 
between $\varprojlim_{f} X(K)$ and  $X(K)_{per,f}$. 

\begin{Lemma}
\label{lemma:elemetary:property:of:periodic:pts}
\begin{enumerate}
\renewcommand{\labelenumi}{(\roman{enumi})}
\item
If $P$ is a $K$-periodic point, 
then there is an element $(x_n)_{n=0}^{\infty} \in \varprojlim_{f} X(K)$ 
such that $P = x_0$. By this correspondence, 
$X(K)_{per, f}$ can be seen as a subset of $\varprojlim_{f} X(K)$.
We say an element of $\varprojlim_{f} X(K)$ which 
lies in the image of $X(K)_{per, f}$ to be periodic.
\item
If $X(K)_{per, f} \subsetneq \varprojlim_{f} X(K)$ 
in the above correspondence, 
then $\varprojlim_{f} X(K)$ is an infinite set.  
\item
If $\varprojlim_{f} X(K)$ is a finite set, 
then $X(K)_{per, f} = \varprojlim_{f} X(K)$ 
in the above correspondence. In particular, 
$X(K)_{per,f}$ is also a finite set.   
\end{enumerate}
\end{Lemma}

\Proof
(i) Let $f^{\circ p}(P) = P$. For any $n \geq 0$, 
we divide $n$ by $p$ and write 
$n = q p + l$ with $0 \leq l \leq p-1$. 
Then if we put $x_n = f^{\circ (p-l)}(P)$,  
$(x_n)_{n=0}^{\infty}$ is an element of $\varprojlim_{f} X(K)$.   

(ii)
Suppose $(x_n)_{n=0}^{\infty} \in \varprojlim_{f} X(K)$ is not periodic.  
By lemma~\ref{lemma:useful}, for any fixed $m$, there are only finitely 
many $k$ with $x_k = x_m$. Then $\{ (x_n)_{n=m}^{\infty} \mid m \geq 0 \} 
\subset \varprojlim_{f} X(K)$ 
is an infinite set. 

(iii)
If $\varprojlim_{f} X(K)$ is a finite set, then 
every $(x_n)_{n=0}^{\infty} \in \varprojlim_{f} X(K)$ is periodic 
by (ii). In particular, $x_0$ is periodic. Therefore, 
the correspondence of (i) becomes bijective. 
\QED

Next lemma shows that finiteness still holds  
if we change $f$ to some powers of $f$. 

\begin{Lemma}
\label{lemma:power}
Let $k$ be a positive integer. 
\begin{enumerate}
\renewcommand{\labelenumi}{(\roman{enumi})}
\item
$X(K)_{per,f^{\circ k}}$ is a finite set 
if and only if  $X(K)_{per,f}$ is a finite set. 
\item
$\varprojlim_{f^{\circ k}} X(K)$ is a finite set
if and only if $\varprojlim_{f} X(K)$ is a finite set. 
\end{enumerate}
\end{Lemma}

\Proof
(i)
Suppose $P$ satisfies $f^{\circ m}(P) = P$. Then 
$P$ satisfies $(f^{\circ k})^{\circ m}(P) = P$. 
This shows that 
$X(K)_{per,f} = X(K)_{per,f^{\circ k}}$. 

(ii)
We have only to prove the ``only if''part. 
If $\varprojlim_{f^{\circ k}} X(K)$ is a finite set, 
its elements are all periodic by 
Lemma~\ref{lemma:elemetary:property:of:periodic:pts}(ii). 
Thus if we set 
\[
S = \{ x \in X(K) \mid 
      \text{there is an 
 $(x_n)_{n=0}^{\infty} \in \varprojlim_{f^{\circ k}} X(K)$
 and an $m$ such that $x = x_m$.} \},
\]
then $S$ is a finite set. Now the finiteness of 
$\varprojlim_{f} X(K)$ follows from Lemma~\ref{lemma:useful}.
\QED

Now we prove the following theorem.

\begin{Theorem}
\label{theorem:finite}
Let $X$ be a projective variety defined over a number field $K$ 
and $f : X \to X$ a surjective 
morphism defined over $K$. 
Assume that there is an ample line bundle 
$L$ such that $f^{*}(L) \otimes L^{-1}$ is ample.  
Then $\varprojlim_{f} X(K)$ is a finite set. 
In particular, $X(K)_{per,f}$ is also a finite set and 
$ \# \varprojlim_{f} X(K) = \# X(K)_{per,f}$.  
\end{Theorem}

\Proof
If we take a positive rational number $\epsilon'$ 
which is sufficiently small, 
then $f^{*}(L) \otimes L^{-(1 + \epsilon')}$ is still 
ample as a $\QQ$-line bundle. Then by 
Theorem~\ref{theorem:positiveness:of:height}(i), 
and by the fact that  $h_{f^*(L)}(P) - h_L(f(P))$ is a bounded function, 
we have a constant $C$ such that 
\[
h_{L}\left(f(P)\right) - (1 + \epsilon') h_L(P) \geq C. 
\]
for all $P \in X(\overline{K})$. 
Let us take an $\epsilon$ with $0 < \epsilon < \epsilon'$. 
Then there is a constant $M$ such that 
if $h_{L}(P) > M$, then 
\[
h_{L}\left(f(P)\right) - (1 + \epsilon) h_L(P) > 0. 
\]

Now let us define a set $S$ to be 
\[
S = \{ x \in X(K) \mid h_{L}(x) \leq M \}. 
\] 
Since $L$ is ample, 
$S$ is a finite set by Northcott.   

In the following we show that, 
if $(x_n)_{n=0}^{\infty} \in \varprojlim_{f} X(K)$, then  
there is a subsequence 
$(x_{n_i})_{i=0}^{\infty}$ consisting of elements in $S$. 
In fact, suppose on the contrary that there is an $m$ such that, 
for any $n \geq m$, $x_n$ does not belong to $S$. 
Since $h_{L}(x_n) >M$ for $n \geq m$, we have 
\[
\cdots < (1 + \epsilon)^2 h_{L}(x_{m+2}) < (1 + \epsilon) h_{L}(x_{m+1}) 
< h_{L}(x_m)
\]
This is a contradiction because 
\[
h_{L}(x_n)  < \frac{1}{(1 + \epsilon)^{n-m}} h_{L}(x_m) \rightarrow 0
\quad\quad
(n \to \infty)
\]

Now by applying Lemma~\ref{lemma:useful},  
we get that $(x_n)_{n=0}^{\infty}$ is periodic 
and uniquely determined by $x_0$. 
We also get that the number of $\varprojlim_{f} X(K)$ 
does not exceed the number if $S$. This proves the first assertion. 
The second assertion follows from 
Lemma~\ref{lemma:elemetary:property:of:periodic:pts}.  
\QED

As a corollary, we obtain the finiteness for a certain class of 
varieties. 

\begin{Corollary}
\label{cor:finite} 
Let $X$ be a geometrically irreducible normal projective variety 
defined over a number field $M$.  
Assume that the Picard number of $X$ is $1$ 
(for example, $X$ is $\PP^n$ 
or a geometrically irreducible normal projective curve ). 
Then $X$ is periodically finite. 
\end{Corollary}

\Proof
Let $K$ be a finite extension field of $M$ 
and $f : X \to X$ be a surjective $K$-morphism of 
$\deg f \geq 2$. 
We take an arbitrary ample line bundle $L$ on $X$. 
Then by our hypothesis, there is a integer $d \geq 2$ 
such that  
$f^{*}(L)$ is numerically equivalent to 
$L^{\otimes d}$. 
In particular, $f^{*}(L) \otimes L^{-1}$ is ample. 
\QED

Let us keep the notation of Theorem~\ref{theorem:finite}. 
Assume here that $f^*(L)$ is linearly equivalent to $L^{\otimes d}$. 
In this case, due to Tate, there exists a unique height function 
$h_{L, f}$ such that $h_{L, f} = h_L + O(1)$ and that 
$h_{L, f}\left( f(P) \right) = d h_{L}(P)$ 
(cf. \cite[Chap 4. Proposotion 1.9]{La}). 
Then for any periodic points with respect to $f$, 
their height must be $0$ 
with respect to $h_{L, f}$. For example:

\begin{Corollary}
\label{cor:m-plication}
Let $K$ be a number field, 
$A$ an Abelian variety defined over $K$ and 
$[m] : A \to A$ the $m$-plication map with
$m \geq 2$. 
Then $\varprojlim_{[m]} A(K)$ is a finite set 
and the number of $\varprojlim_{[m]} A(K)$ does not exceed  
the number of torsion $K$-points. 
\end{Corollary}

\Proof
Extending $K$ if necessary, we may assume that 
there is an ample symmetric line bundle $L$ on $A$.  
Then $f^{*} L \simeq L^{\otimes m^2}$ and  
we can apply the theorem. In this case, if $x$ is a periodic point, 
then $x$ is a torsion point.
\QED

We finish this section by giving examples such that 
$X(K)_{per,f}$ is infinite. 

\begin{Example}
We give an example such that 
$X(K)_{per,f}$ (and thus $\varprojlim_{f} X(K)$) is infinite. 
Let $E$ be an elliptic curve defined over a number field $K$
such that $E(K)$ is an infinite set.   
Let $X$ be $ E \times E$ and 
$f : X \to X$ map $(P,Q)$ to $(P, [2](Q))$. 
Then $f$ is finite of degree $4$   
and the points of the form 
$(P, 0)$ are all periodic points. 
\end{Example}

\begin{Example}
We give an example such that 
$X(K)_{per,f}$ is finite but $\varprojlim_{f} X(K)$ is infinite.  
Let $E$ be an elliptic curve defined over a 
number field $K$ for which $E(K)$ 
contains non-torsion points. 
Let $P_0 \in E(K)$ be a non-torsion point. 
Let $X$ be $ E \times E$ and 
$f : X \to X$ map $(P,Q)$ to $(P+P_0, [2](Q))$. 
Then $f$ is finite of degree $4$ and 
contains a sequence 
$(x_n)_{n=0}^{\infty} \in \varprojlim_{f} X(K)$ 
with $x_n = (-[n](P_0), 0)$. Thus by 
Lemma~\ref{lemma:elemetary:property:of:periodic:pts},
$\varprojlim_{f} X(K)$ is not finite.   
On the other hand, there are no periodic points. 
\end{Example}

We note that we can give examples similar to 
the above two examples by using $\PP^1$. 

\medskip
\section{Curves}
\label{sec:curves}
By a curve, we mean an integral separated scheme 
of finite type over a ground field. In this section, 
we prove that a curve is periodically finite. 
Since there is no surjective morphism $f : C \to C$ 
with $\deg f \geq 2$ if $C$ is a smooth projective curve 
of genus $\geq 2$, we are mainly concerned with a curve $C$ 
such that $C \otimes \overline{\QQ}$ is a reduced scheme 
consisting of rational curves and elliptic curves. 
First we prove two lemmas.

\begin{Lemma}
\label{lemma:general:curve}
Let $C$ be a curve defined over $\overline{\QQ}$, 
and $f : C \to C$ a morphism over $\overline{\QQ}$. 
Then there is a completion $\overline{C}$ of $C$ 
and a morphism $\overline{f} : \overline{C} \to \overline{C}$ 
which is an extension of $f$. 
\end{Lemma}

\proof
Let us take an arbitrary complete curve $\overline{C}'$ 
which is a completion of $C$ 
and set $T' = \overline{C}' \setminus C (\overline{\QQ}$. 
If $t \in T'$ is a singular point of $\overline{C}'$, 
then we blow it up. Iterating this procedure, we get 
a completion $\overline{C}$ such that every point in 
$T = \overline{C} \setminus C (\overline{\QQ})$ is a smooth point of 
$\overline{C}$. Now $f$ defines a rational map 
$\overline{f} : \overline{C} \cdots\to \overline{C}$. 
Since it is defined over $T$ and $C$, 
$\overline{f}$ is actually a morphism. 
\QED

\begin{Lemma}
\label{lemma:smooth:curve}
Let $C$ be a curve defined over a number field $M$ 
which is geometrically irreducible. 
Then $C$ is periodically finite.  
\end{Lemma}

\Proof
Let $K$ be a finite extension of $M$ and 
$f : X \to X$ a surjective morphism defined over $K$ 
with $\deg f \geq 2$.   
By taking a finite extension of $K$ if necessary, 
Lemma~\ref{lemma:general:curve} indicates that 
there is a completion $\overline{C}$ of $C$ 
and a extension $\overline{f}$ of $f$ which are defined over $K$. 
Then $\varprojlim_{f} C(K)$ can be seen as a subset of 
$\varprojlim_{\overline{f}} \overline{C}(K)$. 
For a general point $P \in\overline{C} (\overline{\QQ}$, let 
$L = \OO_{\overline{C}}(P)$. Then, 
since $\deg \overline{f} \geq 2$, $f^*(L) \otimes L^{-1}$ is ample. 
Thus, by Theorem~\ref{theorem:finite}, 
$\varprojlim_{\overline{f}} \overline{C}(K)$ is a 
finite set. 
This proves the lemma. 
\QED

Now we prove the following proposition.

\begin{Proposition}
\label{prop:finiteness:for:complete:scheme}
Let $C$ be a reduced scheme which is 
a chain of geometrically irreducible curves 
over $\overline{\QQ}$. 
Let $f : C \to C$ be a surjective morphism 
such that, for every irreducible component $C_i$ of $C$,  
$f \vert_{C_i} : C_i \to f(C_i)$ has degree $\geq 2$. 
Then for a number field $K \subset \overline{\QQ}$ 
such that $C$ and $f$ are defined over $K$, 
$\varprojlim_{f} C(K)$ is a finite set. 
\end{Proposition}

\Proof
If $K'$ is a extension field of $K$, 
then the finiteness of $\varprojlim_{f} C(K')$ 
implies the finiteness of $\varprojlim_{f} C(K)$. 
Thus to prove the proposition, 
we may take a finite extension of $K$ if necessary. 
Let $C_1, C_2, \ldots, C_l$ be the 
irreducible components of $C$. 
Since $f$ is surjective, the dimension of $f(C_{\alpha})$ 
is $1$ for every $\alpha$. Thus $f$ is seen to 
induce a transposition of the set $C_1, C_2, \ldots, C_l$. 
Then $f^{\circ l!}$ maps $C_{\alpha}$ to $C_{\alpha}$ 
for $1 \leq i \leq l$. 
Let us set $S = \left( \cup_{\alpha \neq \beta} 
C_{\alpha} \cap C_{\beta} \right)_{red}$. 
By Lemma~\ref{lemma:power}, we have only to show that 
$\varprojlim_{f^{\circ l!}} C(K)$ is a finite set. 
We may take a sufficiently large $K$, so that 
$C_{\alpha}$'s and $S$ are all defined over $K$. 
Now let  $(x_n)_{n=0}^{\infty} \in \varprojlim_{f^{\circ l!}} X(K)$. 

{\bf Case 1}
Suppose that there exists a subsequence 
$(x_{n_i})_{i=0}^{\infty}$ consisting of elements in $S$. 
Then by Lemma~\ref{lemma:useful}, 
the number of $(x_n)_{n=0}^{\infty}$  in this case is finite. 

{\bf Case 2}\quad
Suppose that 
there is no subsequence 
$(x_{n_i})_{i=0}^{\infty}$ consisting of elements in $S$. 
Then there is an $\alpha$ such that 
every $x_n$ belongs to $C_{\alpha}$. 
By Lemma~\ref{lemma:general:curve}, 
$\varprojlim_{f^{\circ l!}} C_{\alpha}(K)$ is 
a finite set. Thus the number of 
$(x_n)_{n=0}^{\infty}$  in this case is also finite. 
\QED

As a corollary, we get 

\begin{Corollary}
\label{cor:curve}
Let $C$ be a curve defined over a number field $M$.
Then $C$ is periodically finite. 
\end{Corollary}

\Proof
Let $K$ be a finite extension of $M$ and 
$f : C \to C$ be a surjective $K$-morphism 
with $\deg f \geq 2$.  
Let us consider $C_{\overline{\QQ}}$ and 
let $C_1, C_2, \ldots, C_l$ be its 
irreducible components. By abbreviation, 
$f$ also denotes the induced morphism 
$C_{\overline{\QQ}} \to C_{\overline{\QQ}}$. 
Since $C_1, C_2, \ldots, C_l$ are all conjugate to each other, 
the degree of $f\vert_{C_{\alpha}}$ is greater or equal to $2$
for each $1 \leq \alpha \leq l$. 
Now the assertion follows from 
Proposition~\ref{prop:finiteness:for:complete:scheme}. 
\QED

\medskip
\section{Abelian varieties}
Let $A$ be an abelian variety defined over a number field $M$. 
Recall that $A$ is said to be simple 
if $End(A)_{\QQ}$ is simple. 
In this section, we show that 
$A$ is periodically finite if and only if 
$A$ is simple. 

First we show that if an abelian variety is simple, 
then it is periodically finite. 

\begin{Proposition}
\label{prop:simple:abelian:variety}
Let $A$ a simple abelian variety defined over 
a number field $M$. 
Then $A$ is periodically finite. 
\end{Proposition}

\Proof
Let $K$ be a finite extension field of $M$ and 
$f : X \to X$ a finite $K$-morphism with $\deg f \geq 2$. 
Let us set $B_n = \{ P \in A(K) \mid f^{\circ n}(P) = P \}$. 
We prove the finiteness of $A(K)_{per,f}$ in two steps. 

{\bf Step 1}\quad
We assume here that $f$ is a homomorphism. 
Let us denote by $A(K)_{tor}$ 
the set of $K$-valued torsion points on $A$. 
It is well known that $A(K)_{tor}$ is a finite set 
(cf. Corollary~\ref{cor:m-plication}). 
Since $A$ is simple and $f^{\circ n} \neq 1$, 
$B_n = \Ker(f^{\circ n} - 1)(K)$ is a finite abelian group. 
In particular, $B_n \subset A(K)_{tor}$. 
Thus $A(K)_{per,f} = \cup_{n=1}^{\infty} B_n \subset A(K)_{tor}$ 
is a finite set. 

{\bf Step 2}\quad
Here we treat a general $f$. 
If $B_n = \emptyset$ for $n \geq 1$, then 
we have nothing to prove. 
Thus we assume that there is an $k$ with 
$B_k \neq \emptyset$ and we shall prove 
$A(K)_{per,f}$ is a finite set. 
Since $A(K)_{per,f^{\circ k}} = A(K)_{per,f}$ by 
Lemma~\ref{lemma:power}, 
we may assume that $B_1 \neq \emptyset$. 
We take $x_0 \in B_n$, i.e., $f(x_0) = x_0$. 
We give $A$ another group structure such that 
the identity is $x_0$. We denote this abelian variety by $A'$.  
Since $f$ maps $x_0$ to itself, 
$f$ is a homomorphism of $A'$. 
Therefore, $A'(K)_{per,f}$ is a finite set by Step $1$. 
Since $A$ and $A'$ are identical as a set 
and thus $A(K)_{per,f} = A'(K)_{per,f}$, we are done. 
\QED

Next we show that if $A$ is not simple, then 
$A$ is not periodically finite. 
First we prove the following lemma. 

\begin{Lemma}
\label{lemma:many:K:points:on:abelian:variety}
Let $A$ be an abelian variety defined over a 
number field $M$. Then there exists a 
finite extension field $K$ of $M$ 
such that $A(K)$ is an infinite set. 
\end{Lemma}

\Proof
By Bertini, there is a curve $C$ of genus $\geq 2$ on 
$A_{\overline{\QQ}}$. By Raynaud's theorem \cite{Ray} 
(Manin-Mumford conjecture), $C(\overline{\QQ}) \cap 
A(\overline{\QQ})_{tor}$ is a finite set. 
Take a finite extension field $K$ of $M$ 
such that $C(K)$ contains a non-torsion point $P$. 
Then since $A(K)$ contains $P$, the rank of $A(K)$ 
is positive. 
(The author does not know easier proofs of this 
lemma.)    
\QED

\begin{Proposition}
\label{prop:non:simple:abelian:variety}
Let $A$ be an abelian variety defined over 
a number field $M$. If $A$ is not simple, 
then $A$ is not periodically finite. 
\end{Proposition}

\Proof
Since $A$ is not simple, 
there is an $\overline{\QQ}$-isogeny  $g : B \times C \to X$, 
where $B$ and $C$ are positive-dimensional abelian 
varieties. Let us set $D = \Ker g$, which is a finite 
group of order $d = \# D$.   

We consider a morphism 
\[
[d+1] \times [1] : B \times C \longrightarrow B \times C. 
\]
Since, for a point $(b, c) \in D$, $([d]b, [d]c) = 0$, 
we get $[k+1] \times [1] (b, c) = (b, c)$ for 
any $(b, c) \in D$. 
In particular, $[d+1] \times [1]$ induces a morphism
\[
f : A \longrightarrow A. 
\] 
By the snake lemma, $\Ker([d+1] \times [1]) = \Ker f$, thus 
$f$ is a surjective morphism with $\deg f \geq 2$. 
Now we take a finite extension field $K$ of $M$ such that 
$B$ and $C$ are defined over $K$ and that $C(K)$ 
is an infinite set. Then the infinite set 
\[
g (\{ (0, Q) \in B(K) \times C(K)\})
\]
is contained in $A(K)_{per,f}$. 
\QED

Combining Proposition~\ref{prop:simple:abelian:variety} and 
Proposition~\ref{prop:non:simple:abelian:variety}, 
we get:

\begin{Theorem} 
\label{thm:abelian:variety}
Let $A$ be an abelian variety defined over 
a number field . 
Then $A$ is periodically finite if and only if 
$A$ is simple.
\end{Theorem}

\section{surfaces with non-negative Kodaira dimensions}
\label{sec:non:ruled:surfaces}
In this section we consider 
smooth projective surfaces with non-negative Kodaira dimensions. 

E. Sato and Y. Fujimoto \cite{FS} \cite{FS2} determined 
smooth projective varieties of $\dim = 3$ 
with the non-negative Kodaira dimensions which has a 
non-trivial surjective endomorphism. 

As a test case, they consider the surface case, 
which is as in the following. 

\begin{Theorem}[E. Sato and Y. Fujimoto]
If a smooth projective surface $X$ has a 
surjective endomorphism $f :X \to X$ with 
$\deg f \geq 2$, then $X$ must be minimal and  
is one of the following types; 
\begin{enumerate}
\renewcommand{\labelenumi}{(\roman{enumi})}
\item
$X$ is an abelian surface, 
\item
$X$ is a hyperelliptic surface, or
\item
The Kodaira dimension of $X$ is $1$ and 
$X$ carries an elliptic fibration $\pi: X \to B$ 
whose singular fibers are at most multiple of the type ${}_m I_0$
in the sense of Kodaira, 
where $B$ is a smooth projective curve.
\end{enumerate}
\end{Theorem}

\Proof
For the reader's sake, 
we give a brief sketch of a proof. 
 
Since $X$ has non-negative Kodaira dimension, 
$f : X \to X $ must be \'{e}tale 
(cf. \cite[Theorem~11.7]{Iitaka}).  
Suppose there is an exceptional curve $C$ on $X$.  
Then the equality
\[
f^*(C) \cdot K_X = f^*(C) \cdot f^* K_X = - (\deg f) 
\]
shows that there are at least two exceptional curves on $X$. 
Iterating this procedure, we get a contradiction. 

We note that since $f$ is \'{e}tale, 
$\chi_{top}(X) = (\deg f) \chi_{top}(X)$. 
Then $\deg f \geq 2$ implies $\chi_{top}(X) = 0$. 
In the same way, we get $\chi(\OO_X) = 0$.  

If the Kodaira dimension of $X$ is $2$, 
there are no surjective morphisms $f : X \to X$ with 
$\deg f \geq 2$ (cf. \cite[Proposition~10.10]{Iitaka}). 

If the Kodaira dimension of $X$ is $1$, 
then $\chi_{top}(X) = 0$ indicates 
that $X$ has possibly only multiple singular fibers of type ${}_m I_0$. 

If the Kodaira dimension of $X$ is $0$, 
then $\chi(\OO_X) = 0$ indicates that 
$X$ cannot be a K3 surface nor an Enriques surface. 
\QED

We determined in the previous section 
when an abelian surface is periodically finite. 
Now we study whether a surface of the case (ii) or (iii) 
is periodically finite. 

\begin{Proposition}
\label{prop:bielliptic:surface}
Let $X$ be a hyperelliptic surface defined over a 
number field $M$. Then $X$ is not periodically finite. 
\end{Proposition}

\Proof
Let $E$, $F$ be arbitrary elliptic curves, 
$G$ a group of translations of $E$ which operates on $F$. 
According to the Bagnera-De Franchis list 
(\cite[Liste~VI.20]{Beauville}), 
all the hyperelliptic curves are one of the following types;
\begin{enumerate}
\renewcommand{\labelenumi}{(\roman{enumi})}
\item
$X \cong (E \times F) / G$, $G = \ZZ /2$ 
operating on $F$ by $x \mapsto -x$,
\item
$X \cong (E \times F) / G$, $G = \ZZ /2 \times \ZZ /2$ 
operating on $F$ by $x \mapsto -x$, $x \mapsto x + \epsilon$ 
$(\epsilon \in F_2)$,
\item
$X \cong (E \times F_i) / G$, $G = \ZZ /4$ 
operating on $F$ by $x \mapsto ix$, 
where $F_i = \CC / \ZZ + i \ZZ$, 
\item
$X \cong (E \times F_i) / G$, $G = \ZZ /4$ 
operating on $F$ by $x \mapsto ix$,
\item
$X \cong (E \times F_{\rho}) / G$, $G = \ZZ /3$ 
operating on $F$ by $x \mapsto \rho x$, 
where $\rho = \frac{-1 + \sqrt{-3}}{2}$ and 
$F_{\rho} = \CC / \ZZ + \rho \ZZ$. 
\item
$X \cong (E \times F_{\rho}) / G$, $G = \ZZ /3 \times \ZZ /3$ 
operating on $F$ by $x \mapsto \rho x$, 
$x \mapsto X + \frac{1 -\rho}{3}$ 
\item
$X \cong (E \times F_{\rho}) / G$, $G = \ZZ /6$ 
operating on $F$ by $x \mapsto -\rho x$.
\end{enumerate}
Now we consider the case (i). In this case, 
\[
[3] \times [1] : E \times F \longrightarrow E \times F
\]
induces a surjective morphism 
\[
f : X \to X
\]
with $\deg f \geq 2$. If we take a sufficiently large 
finite extension field $K$ of $M$, 
Then the infinite set $\{ (0, Q) \mid Q \in F(K)\}$ is 
contained in $(E \times F)(K)_{per,[3] \times [1]}$. Thus 
$X(K)_{per,f}$ is also an infinite set. 
The other cases can be treated in similar ways. 
In lieu of $[3] \times [1]$, 
we have only to consider $[g+1] \times [1]$ where 
$g = \# G$. 
\QED

Next we treat a case of an elliptic surface. 
We first consider an elliptic surface such that the genus 
of the base curve is greater or equal to $2$, and then  
one such that the genus of the base curve is $0$ or $1$. 

\begin{Proposition}
\label{prop:Zeifert}
Let $M$ be a number field. Let $X$ be a smooth projective 
surface defined over $M$ with the Kodaira dimension $1$. 
We assume that $X$ carries an elliptic fibration 
$\pi : X \to B$ with at most multiple singular fibers of the 
type ${}_m I_0$ in the sense of Kodaira, 
where $B$ is a smooth projective curve of 
genus $g(B) \geq 2$. 
Then $X$ is periodically finite. 
\end{Proposition}

\Proof
Let $f : X \to X$ be a surjective morphism 
with $\deg f \geq 2$. 
Since $X$ has a unique structure of an elliptic fibration up to 
isomorphisms, there is an automorphism 
$g : B \to B$ with $\pi \circ f = g \circ \pi$. 
Moreover there is a positive integer $k$ such that  
$g^{\circ k}$ is the identity morphism. 
By Lemma~\ref{lemma:power}, we may assume by interchanging 
$f$ with $f^{\circ k}$ that $g$ is the identity morphism.  
Then $f$ is a finite morphism preserving fibers. 
Let $K$ be a sufficiently large number field such that 
$X, B, f, \pi, g$ are all defined over $K$. 

Now let $b$ be a point of $B$ and consider 
$f\vert_{(X_b)_{red}} : (X_b)_{red} \to (X_b)_{red}$. 
Since $f$ is an \'{e}tale morphism 
(cf. \cite[Theorem~11.7]{Iitaka}),  
$f\vert_{(X_b)_{red}} : (X_b)_{red} \to (X_b)_{red}$ 
is a surjective morphism with $\deg f\vert_{(X_b)_{red}} \geq 2$. 

Now we prove the finiteness of $X(K)_{per,f}$ by 
showing the finiteness of $\varprojlim_{f} X(K)$ 
(cf. Lemma~\ref{lemma:elemetary:property:of:periodic:pts}). 
Let $(x_n)_{n=0}^{\infty} \in \varprojlim_{f} X(K)$. 
Since $f$ preserves fibers, $x_n$ are all 
contained in the fiber $X_{\pi(x_0)}$. 
On the other hand, by Mordell-Faltings' theorem, $B(K)$ 
is a finite set. Since $\pi(x_0) \in B(K)$,  
the number of $b$ such that $b = \pi(x_0)$ for some 
$(x_n)_{n=0}^{\infty} \in \varprojlim_{f} X(K)$ is finite.  
Since $\varprojlim_{f\vert{(X_b)_{red}}} 
(X_{b})_{red}(K)$ is a finite set 
for each such $b$ 
by Lemma~\ref{lemma:smooth:curve}, 
we get the assertion. 
\QED

Next we consider an elliptic surface 
such that the genus of the base curve is $0$ or $1$. 
We prove the following two lemmas in advance. 

\begin{Lemma}
\label{lemma:no:small:multi:section}
Let $\pi : X \to B$ be an elliptic surface with the
Kodaira dimension $1$. 
Then the geometric genus of every multi-section of $\pi$ 
is greater or equal to $2$. 
\end{Lemma}

\Proof
Suppose there is a multi-section $C$ on $X$ 
such that the geometric genus of $C$ is $0$ or $1$. 
Then there is an elliptic curve $B'$ with a 
surjection $u : B' \to C$. 
Let us set $v = u \circ \pi : B' \to B$. 
Now we consider the following Cartesian product, 
\[
\begin{CD}
X' @>{v'}>> X \\
@V{\pi'}VV @V{\pi}VV \\
B' @>{v}>> B.
\end{CD}
\]
Since the singular fibers of $\pi'$ are at most 
multiple fibers of type ${}_m I_0$ 
and since $\pi' : X' \to B'$ has a section, 
$\pi'$ must be a smooth morphism. 
Then there is an elliptic curve $B''$ and 
an \'{e}tale covering $B'' \to B'$ such that    
its pull-back $\pi'' : X'' = X' \times_B' B''$ is trivial, 
i.e., $X''$ is a product of elliptic curves
(cf. \cite[Proposition~VI.8]{Beauville}). 
On the other hand, since there is a surjective morphism
$X'' \to X$, the Kodaira dimension of $X''$ must be greater or 
equal to $1$. This is a contradiction.  
\QED

\begin{Lemma}
\label{lemma:periodic:pt:is:torsion:point}
Let $A$ be a simple abelian variety defined over a 
algebraically closed field $k$ and 
$f : A \to A$ a surjective morphism with $\deg f \geq 2$. 
Assume that there is a positive integer $l$ 
such that $f^{\circ l}(0) = 0$. 
Then $f$ maps a torsion point to a torsion point.   
\end{Lemma}

\Proof
Let us set $f(0) = a$. 

We first show that $a$ is a torsion point. 
For this purpose, we may assume that $l \geq 2$. 
Let us set $h = f - a : A \to A$. 
Then $h$ is a homomorphism with $\deg h \geq 2$. 
The equality $f^{\circ l}(0) = 0$ indicates that 
\[
h^{\circ (l-1)}(a) + h^{\circ (l-2)}(a) + \cdots + h(a) + a = 0.
\]
Let us set 
\[
h' = h^{\circ (l-1)} + h^{\circ (l-2)} + \cdots + h + id.
\]

We claim that $h' : A \to A$ is a surjective homomorphism 
of $\deg h' \geq 2$. Indeed, suppose $h'$ is not surjective. 
Then, since $A$ is simple, $h'$ must be the zero map. 
Then we get 
\[
- id = h \circ \left( h^{\circ (l-2)} + h^{\circ (l-2)} + \cdots + id
               \right), 
\] 
which contradicts with $\deg h \geq 2$. 
Thus $h'$ is surjective. Moreover, since $h'$ maps a non-zero element $a$ 
to $0$, $h'$ is not an isomorphism.   

Then, since $A$ is simple, $\Ker h'$ must be a finite group. 
In particular, $a \in \Ker h'$ is a torsion point. 

Now let $b$ be any torsion point of $A$. 
We take a positive integer $n$ such that 
$n a = 0$ and $n b =0$. 
Then 
\begin{align*}
n f(b) & = n (h (b) + a) \\
       & = h (n b) + n a = 0. 
\end{align*}
Thus $f$ maps a torsion point to a torsion point. 
\QED

\begin{Remark}
The above lemma does not hold in general 
for abelian varieties. 
For example, 
let $A$ be an abelian variety and $a$ a non-torsion point. 
If we set $f : A \times A \to A \times A$ by 
$f(x, y) = (2x, -y+a)$, then $f^{\circ 2}(0, 0) = (0,0)$. 
However, $f(0,0) = (0,a)$ is not a torsion point. 
\end{Remark}

\begin{Proposition}
\label{prop:Zeifert:0:1}
Let $M$ be a number field. Let $X$ be a smooth projective 
surface defined over $M$ with the Kodaira dimension $1$. 
We assume that $X$ carries an elliptic fibration 
$\pi : X \to B$ with at most multiple singular fibers of the 
type ${}_m I_0$ in the sense of Kodaira, 
where $B$ is a smooth projective curve of 
genus $0$ or $1$. 
Then $X$ is periodically finite. 
\end{Proposition}

\Proof
Let $f : X \to X$ be a surjective morphism 
with $\deg f \geq 2$. 
Since $X$ has a unique structure of an elliptic fibration up to 
isomorphisms, there is an automorphism 
$g : B \to B$ with $\pi \circ f = g \circ \pi$. 
Let $K$ be a sufficiently large number field such that 
$X, B, f, \pi, g$ are all defined over $K$. 

{\bf Case 1}\quad
Suppose that for any $k \geq 1$, $g^{\circ k}$ is not the identity 
morphism. Let us set 
\[
S = \{ b \in B(\overline{\QQ}) \mid 
       \text{$g^{\circ k} (b) = b$ for some $k \geq 1$}
    \}. 
\]

We claim that $S$ consists at most two points. 
Indeed, suppose $S$ contains three points 
$b_1, b_2, b_3 \in B(\overline{\QQ})$ 
such that $g^{\circ k_i} (b_i) = b_i$ for 
$i = 1, 2, 3$. Then for $k = k_1 k_2 k_3$ we get 
$g^{\circ k} (b_i) = b_i$ for $i = 1, 2, 3$. 
Since $B$ is $\PP^1$ or an elliptic curve, 
this shows that $g^{\circ k}$ is the identity morphism, 
which contradicts our assumption of Case 1. 

We take $l$ such that $g^{\circ l}(b) = b$ for any $b \in S$.  
Now we prove the finiteness of $X(K)_{per,f}$ by 
showing the finiteness of $\varprojlim_{f^{\circ l}} X(K)$ 
(cf. Lemma~\ref{lemma:elemetary:property:of:periodic:pts} and 
Lemma~\ref{lemma:power}). 
Let $(x_n)_{n=0}^{\infty}$ be an element of 
$\varprojlim_{f^{\circ l}} X(K)$. 
Since $\pi(x_n)$ belongs to $S$, $x_n$ are all 
contained in the fiber $X_{\pi(x_0)}$. 
Since $f^{\circ l}$ is an \'{e}tale morphism 
(cf. \cite[Theorem~11.7]{Iitaka}), 
$\varprojlim_{f^{\circ l} \vert{(X_b)_{red}} }
(X_{b})_{red}(K)$ is a finite set for $b \in S$   
by Lemma~\ref{lemma:smooth:curve}. 
Using the finiteness of $S$, we obtain 
the finiteness of $\varprojlim_{f^{\circ l}} X(K)$. 

{\bf Case 2}\quad
Suppose that there is a $k \geq 1$ such that 
$g^{\circ k}$ is the identity morphism. 
To prove the finiteness of $X(K)_{per,f}$, 
we may (and will) assume by interchanging 
$f$ with $f^{\circ k}$ that $g$ is the identity morphism 
(cf.  Lemma~\ref{lemma:power}).  

To prove the theorem, 
we first recall the Merel theorem (cf. \cite{Me}) : 
There is a positive integer $m_0 = m_0(K)$ which depends only
on $K$ such that for any elliptic curve $E$ defined 
over $K$, 
\[
\# E(K)_{tor} \leq m_0(K). 
\]

\begin{Claim}
There is a positive integer $m = m(K)$ which depends only 
on $K$ such that if  $x$ belongs to  $X(K)_{per, f}$,  
then $f^{\circ m}(x) = x$. 
\end{Claim}

\Proof
Let $x$ be an element of $X(K)_{per,f}$. 
Let us set $\pi (x) = b \in B(K)$. 
We introduce a group structure on $X_b$ by 
letting $x$ be the origin. 
Then $X_b$ is an elliptic curve defined over $K$. 
Since $x$ is a periodic point, 
there is a positive integer $l$ 
such that $f^{\circ l} (x) = x$. 
Then by Lemma~\ref{lemma:periodic:pt:is:torsion:point}, 
$f$ maps the set $X_b(K)_{tor}$ to itself. 
On the other hand, by the Merel theorem, 
\[
\# X_b(K)_{tor} \leq m_0. 
\]
Thus, if we set $m = m_0(K) !$, then $m$ depends only on $K$ and 
$f^{\circ m}(x) = x$. 
\QED

Going back to the proof of Proposition~\ref{prop:Zeifert:0:1}, 
we define a reduced subscheme $C$ of $X$ by
\[
C = \{ x \in X \mid f^{\circ m}(x) = x 
    \}.
\]
By the above claim, 
$X(K)_{per, f} \subset C(K)$. 
On the other hand, let 
\[
C = C_1 \cup \cdots \cup C_{\alpha}
\]
be the irreducible decomposition of $C$. 
Since $f$ is \'{e}tale with $\deg f \geq 2$, 
$C$ does not contain a fibral curve. 
If $C_i$ is a horizontal curve, 
then $C_i(K)$ is a finite set by Lemma~\ref{lemma:no:small:multi:section}
and the Mordell-Faltings theorem. 
Therefore $C(K)$ and thus $X(K)_{per, f}$ is a finite set. 
\QED

Combining all the results of this section, we get:

\begin{Theorem}
\label{theorem:surfaces:with:kx:geq:0}
Let $X$ be a smooth projective surface 
with the non-negative Kodaira dimension 
such that $X$ is defined over 
a number field. 
Then $X$ is {\em not} periodically finite if and only if 
$X$ is one of the following types;
\begin{enumerate}
\renewcommand{\labelenumi}{(\roman{enumi})}
\item
$X$ is an abelian surface which is {\em not} simple, or
\item
$X$ is a hyperelliptic surface. 
\end{enumerate}
\end{Theorem}

\medskip
\section{Fields of finite type over $\QQ$}
\label{sec:final}
In this section, we work over a field of 
finite type over $\QQ$. 
A. Moriwaki has recently constructed the theory 
of height functions over a field of finite type over 
$\QQ$. We first recall a part of his theory. 
We refer to \cite{Mo} for details.  

Let $K$ be a field of finite type over 
$\QQ$ with 
$\td_{\QQ}(K) = d$. Let $B$ be a normal variety which is 
projective and flat over $\ZZ$ such that the 
field of rational functions of $B$ is $K$. 
Let $\overline{H} = (H, h_H)$ be 
a nef $C^{\infty}$-hermitian line bundle on $B$, 
i.e., $H$ is a line bundle on $B$ and $h_H$ is a 
$C^{\infty}$-hermitian line bundle such that for any 
curve on $C$ on $B$, 
$\adeg\left(\achern_{1}(\overline{H}\vert_{C})\right) \geq 0$ 
(in the sense of the Arakelov geometry) and 
that the Chern form $\chern_1(\overline{H})$ is semi-positive. 
There exist many such $\overline{B} = (B, \overline{H})$.  
We pick up a $\overline{B}$ and 
fix it in the following. 

Now, for a point $x \in \PP^n(\overline{K})$, 
let us define $h^{\overline{B}}(x)$ to be  
\[
h^{\overline{B}}(x) = \sum_{\Gamma} 
         \log \left( \max_{1 \leq i \leq n} \{ -\ord_{\Gamma}(x_i) \}
         \adeg\left(\achern_{1}(\overline{H}\vert_{\Gamma})^d\right)
              \right) 
+ \int_{B(\CC)} \log \left( \max_{1 \leq i \leq n} 
                 \{ \vert x_i\vert \}
              \right)\chern_{1}(\overline{H})^d , 
\] 
where 
$x = (x_0, x_1, \ldots, x_n) \in \PP^n(K')$ is its coordinate 
over a sufficiently large extension field $K'$ of $K$, 
and $\Gamma$ runs through all prime divisors on $B$. 
This gives rise to a function 
$h^{\overline{B}} : \PP^n(\overline{K}) \to \RR$. 

Now let $X$ be a projective variety defined over $K$, 
$\phi : X \to \PP^n$ a morphism over $K$. 
For a point $x \in X(\overline{K})$, 
we define the height of $x$ with respect to $\phi$, 
denoted by $h^{\overline{B}}_{\phi}(x)$, to be 
$h^{\overline{B}}_{\phi}(x) = h(\phi(x))$. 

Then the following theorem holds as is the number field case 
(cf. \cite[\S 3 - \S 4]{Mo}). 

\begin{Theorem}
For every line bundle $L$ on a projective variety $X$ 
defined over $K$, 
there exists a unique function 
$h^{\overline{B}}_L : X(\overline{K}) \to  \RR$ modulo bounded functions 
with the following property; 
\begin{enumerate}
\renewcommand{\labelenumi}{(\roman{enumi})}
\item 
For any two line bundles $L_1, L_2$, 
$h^{\overline{B}}_{L_1 \otimes L_2} 
= h^{\overline{B}}_{L_1} + h^{\overline{B}}_{L_2} + O(1)$.  
\item 
If $f : X \to Y$ be a morphism of projective varieties 
over $K$, 
then $h^{\overline{B}}_{f^*(L)} = f^*(h^{\overline{B}}_L) + O(1)$.
\item 
If $\phi : X \to \PP^n$ a morphism over $K$, 
then $h^{\overline{B}}_{\phi^*\left(\OO_{\PP^n}(1)\right)} 
= h^{\overline{B}}_{\phi} + O(1)$. 
\end{enumerate}
Moreover the following properties hold. 
\begin{enumerate}
\renewcommand{\labelenumi}{(\alph{enumi})}
\item
(positiveness)
If we denote 
$\Supp \left( \Coker (H^0(X,L) \otimes \OO_X) \to L \right)$ 
by $\Bs (L)$, then 
$h^{\overline{B}}_L$ is bounded below on 
$\left(X \setminus \Bs(L) \right)$. 
\item
(Northcott)
Assume $L$ is ample.  
Then for any $e \geq 1$ and $M \geq 0$, 
\[
\{ x \in X(\overline{K}) 
\mid  h^{\overline{B}}_{L}(x) \leq M, \quad [K(x):K] \leq e \}
\] 
is a finite set. 
\end{enumerate}
\end{Theorem}

Aside from the Northcott finite theorem, we used three big theorems. 

The first one is the Mordell-Faltings theorem
(cf. Proposition~\ref{prop:Zeifert}). 
It is known that this is also true for 
a finitely generated field over $\QQ$ (cf. \cite[Chapter VI]{FW}). 

The next one is the Raynaud theorem
(cf. Lemma~\ref{lemma:many:K:points:on:abelian:variety}).
This is actually proven 
for a field of finite type over $\QQ$

The last one is the Merel theorem 
(cf. Proposition~\ref{prop:Zeifert:0:1}).  
Unfortunately this is not known for a field of finite type over $\QQ$. 
Thus Theorem~\ref{theorem:surfaces:with:kx:geq:0} 
must be replaced by the following weaker theorem 
for a field of finite type over $\QQ$. 

\begin{Theorem}
\label{theorem:surfaces:with:kx:geq:0:for:general:case}
Let $X$ be a smooth projective surface 
with the non-negative Kodaira dimension 
such that $X$ is defined over 
a field of finite type over $\QQ$. 
Assume that $X$ does not carry an elliptic fibration 
$X \to B$ with $g(B) \leq 1$, where $g(B)$ denotes the 
genus of the base curve $B$. 
Then $X$ is {\em not} periodically finite if and only if 
$X$ is one of the following types;
\begin{enumerate}
\renewcommand{\labelenumi}{(\roman{enumi})}
\item
$X$ is an abelian surface which is {\em not} simple, or
\item
$X$ is a hyperelliptic surface. 
\end{enumerate}
\end{Theorem}

Aside from Proposition~\ref{prop:Zeifert:0:1} and 
Theorem~\ref{theorem:surfaces:with:kx:geq:0:for:general:case}, 
all the other results before this section also hold 
for a field of finite type over $\QQ$.

\end{document}

%%% Local Variables: 
%%% mode: latex
%%% TeX-master: t
%%% TeX-master: t
%%% End: 